\documentclass[a4paper,notitlepage,11pt]{article}%
\usepackage{amssymb}
\usepackage{graphicx}
\usepackage{amsmath}
\usepackage{amsthm}
\usepackage{amsfonts}%
\providecommand{\U}[1]{\protect\rule{.1in}{.1in}}
\theoremstyle{plain}
\newtheorem{theorem}{Theorem}[section]

\newtheorem{corollary}[theorem]{Corollary}

\newtheorem{lemma}[theorem]{Lemma}

\theoremstyle{definition}

\newtheorem{remark}[theorem]{Remark}
\numberwithin{equation}{section}
\numberwithin{theorem}{section}
\ifx\pdfoutput\relax\let\pdfoutput=\undefined\fi
\newcount\msipdfoutput
\ifx\pdfoutput\undefined\else
\ifcase\pdfoutput\else
\ifx\paperwidth\undefined\else
\ifdim\paperheight=0pt\relax\else\pdfpageheight\paperheight\fi
\ifdim\paperwidth=0pt\relax\else\pdfpagewidth\paperwidth\fi
\fi\fi\fi
\begin{document}

\title{Positive solutions of indefinite semipositone problems via sub-super solutions
\thanks{2010 \textit{Mathematics Subject Clasification}. 35J25, 35J60, 35B09.}
 }
\author{U. Kaufmann\thanks{FaMAF, Universidad Nacional de C\'{o}rdoba, (5000)
C\'{o}rdoba, Argentina. \textit{E-mail address: }kaufmann@mate.uncor.edu} , H.
Ramos Quoirin\thanks{Universidad de Santiago de Chile, Casilla 307, Correo 2,
Santiago, Chile. \textit{E-mail address: }humberto.ramos@usach.cl
(Corresponding Author)}
\and \noindent}
\maketitle

\begin{abstract}
Let $\Omega\subset\mathbb{R}^{N}$, $N\geq1$, be a smooth bounded domain, and
let $m:\Omega\rightarrow\mathbb{R}$ be a possibly sign-changing function. We
investigate the existence of positive solutions for the semipositone problem
\[
\left\{
\begin{array}
[c]{lll}%
-\Delta u=\lambda m(x)(f(u)-k) & \mathrm{in} & \Omega,\\
u=0 & \mathrm{on} & \partial\Omega,
\end{array}
\right.
\]
where $\lambda,k>0$ and $f$ is either sublinear at infinity with $f(0)=0$, or
$f$ has a singularity at $0$. We prove the existence of a positive solution
for certain ranges of $\lambda$ provided that the negative part of $m$ is
suitably small. Our main tool is the sub-supersolutions method, combined with
some rescaling properties.

\end{abstract}

\section{Introduction}

Let $\Omega$ be a bounded and smooth domain of $\mathbb{R}^{N}$ with $N\geq1$.
We are concerned with the problem
\[
\left\{
\begin{array}
[c]{lll}%
-\Delta u=\lambda m(x)(f(u)-k) & \mathrm{in} & \Omega,\\
u=0 & \mathrm{on} & \partial\Omega,
\end{array}
\right.  \leqno{(P_\lambda)}
\]
where $\lambda,k>0$ and $m:\Omega\rightarrow\mathbb{R}$ may change sign. Here
either $f$ is nonnegative and continuous on $[0,\infty)$ with $f(0)=0$ or $f$
is positive and continuous in $(0,\infty)$ and singular at $0$. Our main
purpose is to establish existence results for \textit{strictly} positive
solutions of $(P_{\lambda})$.

In the first case, we assume that $f$ is sublinear at infinity. More
precisely, we suppose that $f$ satisfies one of the following conditions:

\begin{itemize}
\item[$(F1)$] For some $0<p<1$,
\begin{equation}
\lim_{s\rightarrow\infty}\frac{f(s)}{s^{p}}=1. \label{sub}%
\end{equation}

\item[$(F2)$] It holds that%
\begin{equation}
\lim_{s\rightarrow\infty}f\left(  s\right)  :=c<\infty. \label{tado}%
\end{equation}

\end{itemize}

Note that $(P_{\lambda})$ belongs in this case to the class of
\textit{semipositone} problems (for applications and further references on
semipositone problems we refer the reader to \cite{costa, dancer} and the nice
survey papers \cite{survey, lee} and references therein). We observe that
proving the existence of a positive solution of $(P_{\lambda})$ is not
straightforward, since the strong maximum principle does not apply. In fact,
even the existence of nonnegative solutions via standard variational methods
is not clear for $(P_{\lambda})$ (even if $m\geq0$).

In \cite{cgs}, the problem $(P_{\lambda})$ was investigated when $m\equiv1$
and $\displaystyle\lim_{s\rightarrow\infty}\frac{f(s)}{s}=0$. Using the
sub-supersolutions method, the existence of a nonnegative solution for
$\lambda$ sufficiently large was proved. Still in the case $m\equiv1$, the
authors of \cite{costa} followed a non-standard variational approach and
proved non-existence, existence, and multiplicity results for nonnegative
solutions in both the sublinear and superlinear cases. In addition, in the
sublinear case they proved the existence of a positive solution when $\Omega$
is a ball.

On the other hand, in \cite{ambro} the authors employed bifurcation and degree
theory arguments to obtain a positive solution of $(P_{\lambda})$ for
$\lambda$ sufficiently large. However, their result holds when $m>0$ in
$\Omega$. More concretely, they considered the problem
\begin{equation}
\left\{
\begin{array}
[c]{lll}%
-\Delta u=\lambda f(x,u) & \mathrm{in} & \Omega,\\
u=0 & \mathrm{on} & \partial\Omega,
\end{array}
\right.  \label{eambro}%
\end{equation}
and proved the following:

\begin{theorem}
\cite{ambro} \label{ambro} Let $f:\overline{\Omega}\times\lbrack
0,\infty)\rightarrow\mathbb{R}$ be a continuous function such that:

\begin{enumerate}
\item $f(x,0)<0$ for all $x \in\Omega$.

\item There exists a continuous function $b:\overline{\Omega}\rightarrow
(0,\infty)$ such that
\[
\lim_{s\rightarrow\infty}\frac{f(x,s)}{s^{q}}=b(x)\quad\text{uniformly in
}x\in\overline{\Omega},\text{ with }0\leq q<1.
\]

\end{enumerate}

Then there exists $\lambda_{0}>0$ such that \eqref{eambro} has a positive
solution for $\lambda\geq\lambda_{0}$. More precisely, there exists a
connected set of positive solutions of \eqref{eambro} bifurcating from
infinity at $\lambda_{\infty}=\infty$.
\end{theorem}

Let us mention that they also dealt with the superlinear case, obtaining now a
positive solution for all $\lambda>0$ sufficiently small. Still with $f$
superlinear, but now with $m$ allowed to change sign, one may combine
rescaling and continuity arguments to show that $(P_{\lambda})$ has a positive
solution for $\lambda$ sufficiently small, see \cite{humberto}. These
arguments rely on the fact that the strong maximum principle applies to the
problem
\[
\left\{
\begin{array}
[c]{lll}%
-\Delta u=m(x)u^{p} & \mathrm{in} & \Omega,\\
u=0 & \mathrm{on} & \partial\Omega,
\end{array}
\right.
\]
with $p>1$. However, when $0<p<1$ and $m$ changes sign in $\Omega$, this is no
longer true. We shall then follow a different approach, based on the
sub-supersolutions method.

Let us denote by $\mathcal{P}^{\circ}$ the interior of the positive cone in
$\mathcal{C}_{0}^{1}(\overline{\Omega})$, i.e.,
\[
\mathcal{P}^{\circ}:=\{u\in\mathcal{C}^{1}(\overline{\Omega}):u>0\text{ in
}\Omega\text{, }u=0\text{ on }\partial\Omega\text{ and }\partial u/\partial
\nu<0\text{ on }\partial\Omega\},
\]
where $\nu$ is the outward unit normal to $\partial\Omega$. We also write as
usual $m=m^{+}-m^{-}$ with $m^{+}:=\max\left(  m,0\right)  $ and $m^{-}%
:=\max\left(  -m,0\right)  $.

We next state our first main result:

\begin{theorem}
\label{t1} Let $m\in L^{\infty}\left(  \Omega\right)  $ and $f:[0,\infty
)\rightarrow\lbrack0,\infty)$ be a continuous function with $f(0)=0$. Assume
one of the following conditions:

\begin{enumerate}
\item $f$ satisfies $(F1)$.

\item $f$ satisfies $(F2)$ with $c>k$.
\end{enumerate}

If $\Vert m^{-}\Vert_{L^{\infty}\left(  \Omega\right)  }$ is sufficiently
small, then there exists $\lambda_{0}>0$ such that $(P_{\lambda})$ has a
solution $u_{\lambda}\in\mathcal{P}^{\circ}$ for all $\lambda\geq\lambda_{0}$.
Moreover, $u_{\lambda}(x)\rightarrow\infty$ as $\lambda\rightarrow\infty$ for
every $x\in\Omega$.
\end{theorem}

To be more specific, the smallness condition on $m^{-}$ in Theorem \ref{t1}
should hold with respect to $m^{+}$. In case (1), it must be such that, for
some $\delta>0$, the problem
\begin{equation}
\left\{
\begin{array}
[c]{lll}%
-\Delta w=\left(  1-\delta\right)  m^{+}\left(  x\right)  w^{p}-m^{-}\left(
x\right)  w^{p} & \mathrm{in} & \Omega,\\
w=0 & \mathrm{on} & \partial\Omega,
\end{array}
\right.  \label{U}%
\end{equation}
has a solution $w\in\mathcal{P}^{\circ}$, see Corollary \ref{c1} below. Some
examples of such condition on $m$ can be found in \cite[Theorems 2.1 (ii), 3.1
and 3.3]{nodea} when $N=1$ or when $\Omega$ is a ball and $m$ is radial, and
more generally for any smooth bounded domain in \cite[Theorem 3.1]{ans}. We
note that the solution found in the above theorem in fact lies in
$W^{2,q}\left(  \Omega\right)  $ for some $q>N$, and thus also in
$\mathcal{C}^{1,\theta}(\overline{\Omega})$ with $0<\theta<1$.

In case (2), we require $m^{-}$ to be sufficiently small so that the unique
solution $u$ of
\begin{equation}
\left\{
\begin{array}
[c]{lll}%
-\Delta u=m(x) & \mathrm{in} & \Omega,\\
u=0 & \mathrm{on} & \partial\Omega,
\end{array}
\right.  \label{lineal}%
\end{equation}
belongs to $\mathcal{P}^{\circ}$ (regarding this issue, see Remark \ref{eme}
below). Let us point out that the above condition was used in \cite{dai-} in
order tu study existence of positive solutions to some superlinear elliptic
problems, and also \cite{jesusultimo} for getting positive solutions for some
indefinite singular sublinear problems.

Let us also remark that the main difficulty in the proof of Theorem \ref{t1}
is to provide a positive subsolution of $(P_{\lambda})$. Even though Theorem
\ref{t1} is already known in the case $m>0$, to the best of our knowledge, our
way of getting a positive subsolution is new even in this case, and extends
naturally to the case where $m$ changes sign.

Next, we deal with $(P_{\lambda})$ in the singular case, namely, when
$f:(0,\infty)\rightarrow(0,\infty)$ is continuous and satisfies:

\begin{itemize}
\item[$(F3)$] For some $p>0$,
\begin{equation}
\lim_{s\rightarrow0^{+}}f(s)s^{p}=1. \label{sing}%
\end{equation}
In this case we have:
\end{itemize}

\begin{theorem}
\label{t2} Let $m\in L^{\infty}\left(  \Omega\right)  $ and $f:(0,\infty
)\rightarrow(0,\infty)$ be a continuous function satisfying $(F3)$. If $\Vert
m^{-}\Vert_{L^{\infty}\left(  \Omega\right)  }$ is sufficiently small, then
there exists $\lambda_{0}>0$ such that $(P_{\lambda})$ has a positive solution
$u_{\lambda}\in\mathcal{C}^{1}(\Omega)\cap\mathcal{C}(\overline{\Omega})$ for
all $0<\lambda\leq\lambda_{0}$. Moreover, $\displaystyle$%
\[
\lim_{\lambda\rightarrow0^{+}}\Vert u_{\lambda}\Vert_{L^{\infty}\left(
\Omega\right)  }=0.
\]

\end{theorem}

The smallness condition on $m^{-}$ in Theorem \ref{t2} must be such that the
problem
\[
\left\{
\begin{array}
[c]{lll}%
-\Delta w=\left(  m(x)-\delta\right)  w^{-p} & \mathrm{in} & \Omega,\\
w=0 & \mathrm{on} & \partial\Omega,
\end{array}
\right.
\]
admits a positive solution for some $\delta>0$. Such conditions can be found
in \cite[Corollaries 3.7 and 4.6]{singular}. Let us mention that singular
problems similar to $(P_{\lambda})$, with a nonnegative weight $m$, have been
widely studied in the literature (see e.g. \cite{lee,perera,pepe,zhao} and
references therein), while recently the case $k=0$ and $m$ sign-changing was
considered in \cite{singular}.

\section{Proof of the main results}

Given $m\in L^{q}(\Omega)$, $q>N$, let $u\in W^{2,q}(\Omega)\cap W_{0}%
^{1,q}(\Omega)$ be the unique solution of (\ref{lineal}) and let
$\mathcal{S}:L^{q}(\Omega)\rightarrow W^{2,q}(\Omega)$ be the corresponding
solution operator, i.e., $\mathcal{S}\left(  m\right)  :=u$.

\begin{lemma}
\label{lema}Let $m,h\in L^{q}\left(  \Omega\right)  $, $q>N$, with
$0\not \equiv h\geq0$ in $\Omega$. Suppose that for some $\delta>0$ the
problem \eqref{U} has a solution $w\in\mathcal{P}^{\circ}$. Then there exists
$\beta_{0}>0$ such that for all $\beta\in\left(  0,\beta_{0}\right]  $ there
exists $u_{\beta}\in\mathcal{P}^{\circ}$ solution of
\begin{equation}
\left\{
\begin{array}
[c]{lll}%
-\Delta u=m\left(  x\right)  u^{p}-\beta h\left(  x\right)   & \mathrm{in} &
\Omega,\\
u=0 & \mathrm{on} & \partial\Omega.
\end{array}
\right.  \label{cua}%
\end{equation}

\end{lemma}

\textit{Proof}. First we observe that if $\varphi:=\mathcal{S}\left(
m^{+}\right)  \in\mathcal{P}^{\circ}$ and $t\geq\left\Vert \varphi\right\Vert
_{\infty}^{\frac{p}{1-p}}$, then $t\varphi$ is a supersolution of (\ref{cua}).
Indeed,%
\begin{align*}
-\Delta\left(  t\varphi\right)   &  \geq\left(  t\left\Vert \varphi\right\Vert
_{\infty}\right)  ^{p}m^{+}\left(  x\right)  \\
&  \geq m\left(  x\right)  \left(  t\varphi\right)  ^{p}\geq m\left(
x\right)  \left(  t\varphi\right)  ^{p}-\beta h\left(  x\right)
\qquad\text{in }\Omega.
\end{align*}

On the other hand, let $v:=\mathcal{S}\left(  h\right)  $ and fix $\delta>0$
and $w\in\mathcal{P}^{\circ}$ a solution of \eqref{U}. Then there exists
$\beta_{0}>0$ such that for all $\beta\in\left(  0,\beta_{0}\right]  $ we have
that $\beta v\leq\delta w$ in $\Omega$, and consequently
\begin{equation}
\left(  1-\delta\right)  w\leq w-\beta v\qquad\text{in }\Omega.\label{be}%
\end{equation}
Now, for such $\beta$, employing (\ref{be}) we derive that%
\begin{align*}
-\Delta\left(  w-\beta v\right)   &  =\left(  1-\delta\right)  m^{+}\left(
x\right)  w^{p}-m^{-}\left(  x\right)  w^{p}-\beta h\left(  x\right)  \\
&  \leq\left(  1-\delta\right)  ^{p}m^{+}\left(  x\right)  w^{p}-m^{-}\left(
x\right)  \left(  w-\beta v\right)  ^{p}-\beta h\left(  x\right)  \\
&  \leq m^{+}\left(  x\right)  \left(  w-\beta v\right)  ^{p}-m^{-}\left(
x\right)  \left(  w-\beta v\right)  ^{p}-\beta h\left(  x\right)
\qquad\text{in }\Omega.
\end{align*}
In other words, $w-\beta v$ is a subsolution of (\ref{cua}). Thus, applying
the well-known sub-supersolutions method in the presence of well-ordered weak
sub and supersolutions (see e.g. \cite[Theorem 4.9]{du}) we obtain a solution
$u\in H_{0}^{1}\left(  \Omega\right)  \cap L^{\infty}\left(  \Omega\right)  $.
Moreover, by standard regularity arguments we deduce that $u\in W^{2,q}\left(
\Omega\right)  $, $q>N$, and so, since $w-\beta v\in\mathcal{P}^{\circ}$, we
may conclude also that $u\in\mathcal{P}^{\circ}$. $\blacksquare$

\begin{theorem}
\label{t3} Let $m\in L^{\infty}\left(  \Omega\right)  $ and $f:[0,\infty
)\rightarrow\lbrack0,\infty)$ be a continuous function satisfying $(F1)$ and
$f(0)=0$. Assume in addition that for some $\beta>0$ the problem
\[
\left\{
\begin{array}
[c]{lll}%
-\Delta w=m(x)w^{p}-\beta & \mathrm{in} & \Omega,\\
w=0 & \mathrm{on} & \partial\Omega,
\end{array}
\right.  \leqno{(P_{m,\beta})}
\]
has a positive solution $w\in\mathcal{C}(\overline{\Omega})$. Then there
exists $\lambda_{0}>0$ such that $(P_{\lambda})$ has a positive solution
$u_{\lambda}\in W^{2,q}\left(  \Omega\right)  $, $q>N$, for all $\lambda
\geq\lambda_{0}$. Moreover, $u_{\lambda}(x)\rightarrow\infty$ as
$\lambda\rightarrow\infty$ for every $x\in\Omega$.
\end{theorem}

\textit{Proof}. We set $\underline{u}_{\lambda}:=\lambda^{\frac{1}{1-p}}w$, so
that
\[
-\Delta\underline{u}_{\lambda}=\lambda^{\frac{1}{1-p}}\left(  m(x)w^{p}%
-\beta\right)  .
\]
Thus $\underline{u}_{\lambda}$ is a subsolution of $(P_{\lambda})$ if and only
if
\begin{equation}
\lambda^{\frac{1}{1-p}}\left(  m(x)w^{p}-\beta\right)  \leq\lambda m(x)\left(
f(\lambda^{\frac{1}{1-p}}w)-k\right)  . \label{1}%
\end{equation}
This inequality is equivalent to
\begin{equation}
\frac{\beta}{w^{p}}\geq m\left(  x\right)  \left(  1-\frac{f(\lambda^{\frac
{1}{1-p}}w)}{\left(  \lambda^{\frac{1}{1-p}}w\right)  ^{p}}\right)
+\frac{km\left(  x\right)  }{\left(  \lambda^{\frac{1}{1-p}}w\right)  ^{p}}.
\label{2}%
\end{equation}
Let us set $\varepsilon:=\displaystyle\inf_{x\in\Omega}\frac{\beta}{w(x)^{p}}%
$. Since $\displaystyle\lim_{s\rightarrow\infty}\frac{f(s)}{s^{p}}=1$, there
exists $s_{0}>0$ such that
\[
\left\vert 1-\frac{f(s)}{s^{p}}\right\vert <\frac{\varepsilon}{2\Vert
m\Vert_{\infty}}\quad\text{and}\quad\frac{k}{s^{p}}<\frac{\varepsilon}{2\Vert
m\Vert_{\infty}}\quad\text{for }s>s_{0}.
\]
Hence, for $x\in\Omega$ such that $\lambda^{\frac{1}{1-p}}w(x)>s_{0}$ we have
\[
m\left(  x\right)  \left(  1-\frac{f(\lambda^{\frac{1}{1-p}}w(x))}{\left(
\lambda^{\frac{1}{1-p}}w(x)\right)  ^{p}}\right)  +\frac{km\left(  x\right)
}{\left(  \lambda^{\frac{1}{1-p}}w(x)\right)  ^{p}}<\varepsilon\leq\frac
{\beta}{w(x)^{p}},
\]
so that \eqref{1} holds for such $x$.

Now we set $S:=\displaystyle\sup_{0\leq s\leq s_{0}}|s^{p}-f(s)|$ and fix
$\lambda_{0}>0$ such that
\[
\lambda^{\frac{p}{1-p}}\beta-km(x)>S\Vert m\Vert_{\infty}%
\]
for every $x\in\Omega$ if $\lambda\geq\lambda_{0}$. Consequently, if
$\lambda\geq\lambda_{0}$ and $x\in\Omega$ is such that $\lambda^{\frac{1}%
{1-p}}w(x)\leq s_{0}$, then
\[
\left\vert m\left(  x\right)  \left(  \left(  \lambda^{\frac{1}{1-p}%
}w(x)\right)  ^{p}-f(\lambda^{\frac{1}{1-p}}w(x))\right)  \right\vert \leq
S\Vert m\Vert_{\infty}<\lambda^{\frac{p}{1-p}}\beta-km(x),
\]
which yields \eqref{1}. Therefore, $\underline{u}_{\lambda}$ is a subsolution
of $(P_{\lambda})$ for $\lambda\geq\lambda_{0}$.

On the other hand, let $e:=\mathcal{S}\left(  1\right)  $. We define
$\overline{u}_{\lambda}:=t(e+1)$ with $t>0$. Then $\overline{u}_{\lambda}$ is
a supersolution of $(P_{\lambda})$ if and only if
\[
t\geq\lambda m\left(  x\right)  \left(  f(t(e+1))-k\right)  \quad\text{in
}\Omega,
\]
i.e.
\begin{equation}
1\geq\lambda m\left(  x\right)  \frac{f(t(e+1))-k}{t(e+1)}(e+1)\quad\text{in
}\Omega. \label{3}%
\end{equation}
Since $\displaystyle \lim_{s\rightarrow\infty}\frac{f(s)}{s}=0$, it follows
that given $\varepsilon>0$ there exists $s_{1}>0$ such that $\left\vert
\frac{f(s)-k}{s}\right\vert <\varepsilon$ for every $s>s_{1}$. In particular,
if $t>s_{1}$ then $t(e(x)+1)>t>s_{1}$, so that
\[
\left\vert \frac{f(t(e+1))-k}{t(e+1)}\right\vert <\varepsilon\quad\text{in
}\Omega.
\]
Thus we see that \eqref{3} holds for all $t$ sufficiently large. Taking $t$
larger if necessary, we have $\lambda^{\frac{1}{1-p}}w\leq t(e+1)$ in $\Omega
$, i.e. $\underline{u}_{\lambda}\leq\overline{u}_{\lambda}$ in $\Omega$. We
obtain then, for $\lambda\geq\lambda_{0}$, a solution $u_{\lambda}$ of
$(P_{\lambda})$ satisfying $\underline{u}_{\lambda}\leq u_{\lambda}%
\leq\overline{u}_{\lambda}$. In particular, for every $x\in\Omega$ we have
$u_{\lambda}(x)\geq\lambda^{\frac{1}{1-p}}w(x)\rightarrow\infty$ as
$\lambda\rightarrow\infty$.

The proof is now complete. $\blacksquare$

\qquad

As an immediate consequence of Theorem \ref{t3} and Lemma \ref{lema} we have
the following:

\begin{corollary}
\label{c1} Let $m\in L^{\infty}\left(  \Omega\right)  $ and $f:[0,\infty
)\rightarrow\lbrack0,\infty)$ be a continuous function satisfying $(F1)$ and
$f(0)=0$. Assume in addition that for some $\delta>0$ the problem \eqref{U}
has a solution $w\in\mathcal{P}^{\circ}$. Then there exists $\lambda_{0}>0$
such that $(P_{\lambda})$ has a solution $u_{\lambda}\in\mathcal{P}^{\circ}$
for all $\lambda\geq\lambda_{0}$. Moreover, $u_{\lambda}(x)\rightarrow\infty$
as $\lambda\rightarrow\infty$ for every $x\in\Omega$.
\end{corollary}

We now consider the case where $f$ is bounded:

\begin{theorem}
\label{aco}Let $m\in L^{q}\left(  \Omega\right)  $, $q>N$, and $f:[0,\infty
)\rightarrow\lbrack0,\infty)$ be a continuous function satisfying $(F2)$ with
$c>k$ and $f(0)=0$. Assume in addition that $\mathcal{S}\left(  m\right)
\in\mathcal{P}^{\circ}$. Then there exists $\lambda_{0}>0$ such that
$(P_{\lambda})$ has a solution $u_{\lambda}\in\mathcal{P}^{\circ}$ for all
$\lambda\geq\lambda_{0}$. Moreover, $u_{\lambda}(x)\rightarrow\infty$ as
$\lambda\rightarrow\infty$ for every $x\in\Omega$.
\end{theorem}

\textit{Proof}. Let $c$ be given by $\left(  F2\right)  $ and define
$C:=\displaystyle \sup_{s>0}f\left(  s\right)  $. Observe that $C\in\left[
c,\infty\right)  $. We also set, for $\delta>0$,
\[
\Omega_{\delta}:=\left\{  x\in\Omega:dist\left(  x,\partial\Omega\right)
\geq\delta\right\}
\]
and
\[
M_{\delta}\left(  x\right)  :=\left(  c-k\right)  m\left(  x\right)
-\delta\left\vert m\left(  x\right)  \right\vert \chi_{\Omega_{\delta}%
}-\left(  cm\left(  x\right)  +Cm^{-}\left(  x\right)  \right)  \chi
_{\Omega\setminus\Omega_{\delta}},
\]
where $\chi_{A}$ denotes the characteristic function of $A$.

Since $c>k$, $\mathcal{S}\left(  m\right)  \in\mathcal{P}^{\circ}$ and the
solution operator $\mathcal{S}:L^{q}(\Omega)\rightarrow\mathcal{C}%
^{1}(\overline{\Omega})$ is continuous, there exists $\delta>0$ such that the
problem%
\[
\left\{
\begin{array}
[c]{lll}%
-\Delta w=M_{\delta}\left(  x\right)   & \mathrm{in} & \Omega,\\
w=0 & \mathrm{on} & \partial\Omega,
\end{array}
\right.
\]
admits a solution $w\in\mathcal{P}^{\circ}$. Now, $w\left(  x\right)  \geq
\eta>0$ for all $x\in\Omega_{\delta}$ and some $\eta>0$, and hence taking into
account (\ref{tado}) we see that there exists $\lambda_{0}>0$ such that for
all $\lambda\geq\lambda_{0}$ it holds that%
\begin{align*}
\lambda m\left(  x\right)  \left(  f\left(  \lambda w\right)  -k\right)   &
=\lambda m^{+}\left(  x\right)  \left(  f\left(  \lambda w\right)  -k\right)
-\lambda m^{-}\left(  x\right)  \left(  f\left(  \lambda w\right)  -k\right)
\\
&  \geq\lambda m^{+}\left(  x\right)  \left(  c-k-\delta\right)  -\lambda
m^{-}\left(  x\right)  \left(  c-k+\delta\right)  \\
&  =\lambda M_{\delta}\left(  x\right)  \qquad\text{in }\Omega_{\delta}.
\end{align*}
On the other side, for all $\lambda>0$, we also have that%
\begin{align*}
\lambda m\left(  x\right)  \left(  f\left(  \lambda w\right)  -k\right)   &
=\lambda m^{+}\left(  x\right)  \left(  f\left(  \lambda w\right)  -k\right)
-\lambda m^{-}\left(  x\right)  \left(  f\left(  \lambda w\right)  -k\right)
\\
&  \geq-\lambda m^{+}\left(  x\right)  k-\lambda m^{-}\left(  x\right)
\left(  C-k\right)  \\
&  =\lambda M_{\delta}\left(  x\right)  \qquad\text{in }\Omega\setminus
\Omega_{\delta}\text{.}%
\end{align*}
So, $-\Delta\left(  \lambda w\right)  \leq\lambda m\left(  x\right)  \left(
f\left(  \lambda w\right)  -k\right)  $ in $\Omega$ and therefore $\lambda w$
is a subsolution of $(P_{\lambda})$.

Let now $\varphi:=\mathcal{S}\left(  \left\vert m\right\vert \right)
\in\mathcal{P}^{\circ}$ and $C$ be as above. Given $\lambda\geq\lambda_{0}$,
we choose $t_{0}\geq\lambda C$. Since $C>k$, we get, for every $t\geq t_{0}$,
\begin{align*}
-\Delta\left(  t\varphi\right)   &  \geq\lambda m^{+}\left(  x\right)
C+\lambda m^{-}\left(  x\right)  k\\
&  \geq\lambda m^{+}\left(  x\right)  f\left(  t\varphi\right)  +\lambda
m^{-}\left(  x\right)  k\\
&  \geq\lambda m\left(  x\right)  \left(  f\left(  t\varphi\right)  -k\right)
\qquad\text{in }\Omega,
\end{align*}
and thus $t\varphi$ is a supersolution of $(P_{\lambda})$. Enlarging $t$ if
necessary, so that $t\varphi\geq\lambda w$ in $\Omega$, we obtain then, for
$\lambda\geq\lambda_{0}$, a solution $u_{\lambda}$ of $(P_{\lambda})$
satisfying that $\lambda w\leq u_{\lambda}\leq t\varphi$ in $\Omega$. In
particular, $u_{\lambda}(x)\rightarrow\infty$ as $\lambda\rightarrow\infty$
for every $x\in\Omega$, and this ends the proof of the theorem. $\blacksquare$

\begin{remark}
\label{eme}Let $\delta_{\Omega}\left(  x\right)  :=dist\left(  x,\partial
\Omega\right)  $. If $\Omega:=\left(  a,b\right)  \subset\mathbb{R}$, it
follows directly from Lemma 2.1 in \cite{singular} that if
\[
\max\left\{  \int_{a}^{b}\left(  t-a\right)  m^{-}\left(  t\right)
dt,\int_{a}^{b}\left(  b-t\right)  m^{-}\left(  t\right)  dt\right\}
<\int_{a}^{b}\delta_{\Omega}\left(  t\right)  m^{+}\left(  t\right)  dt,
\]
then $\mathcal{S}\left(  m\right)  \in\mathcal{P}^{\circ}$. For a smooth
bounded domain $\Omega\subset\mathbb{R}^{N}$, Lemma 2.2 in \cite{singular}
allows us to reach the same conclusion provided that
\[
c_{\Omega}\left\Vert m^{-}\right\Vert _{L^{q}\left(  \Omega\right)  }%
<\int_{\Omega}\delta_{\Omega}\left(  x\right)  m^{+}\left(  x\right)  dx,
\]
where $c_{\Omega}>0$ is some constant depending only on $\Omega$.
$\blacksquare$
\end{remark}

\begin{remark}
Let $f$ be merely bounded, define $l_{0}:=\underline{\lim}_{s\rightarrow
\infty}f\left(  s\right)  $ and $l_{1}:=\overline{\lim}_{s\rightarrow\infty
}f\left(  s\right)  $, and suppose that $l_{0}>k$. Then, reasoning as in
Theorem \ref{aco} one can show the existence of positive solutions of
$(P_{\lambda})$ for $\lambda$ large enough provided that $\mathcal{S}\left(
l_{0}m^{+}-l_{1}m^{-}\right)  \in\mathcal{P}^{\circ}$. It is an interesting
open question to see if this remains true just under the hypothesis of Theorem
\ref{aco}, namely, $\mathcal{S}\left(  m\right)  \in\mathcal{P}^{\circ}$.
$\blacksquare$
\end{remark}

We finally study the case of a singular nonlinearity. For the precise
definition of sub and supersolutions in this case we refer to \cite[Section
4]{loc}.

\begin{theorem}
Let $m\in L^{\infty}\left(  \Omega\right)  $ and $f:(0,\infty)\rightarrow
(0,\infty)$ be a continuous function satisfying $(F3)$. Assume that for some
$\delta>0$ the problem
\[
\left\{
\begin{array}
[c]{lll}%
-\Delta w=\left(  m(x)-\delta\right)  w^{-p} & \mathrm{in} & \Omega,\\
w=0 & \mathrm{on} & \partial\Omega,
\end{array}
\right.  \leqno{(P_{m,\delta})}
\]
has a positive solution $w\in\mathcal{C}^{1}(\Omega)\cap\mathcal{C}%
(\overline{\Omega})$. Then there exists $\lambda_{0}>0$ such that
$(P_{\lambda})$ has a positive solution $u_{\lambda}\in\mathcal{C}^{1}%
(\Omega)\cap\mathcal{C}(\overline{\Omega})$ for every $0<\lambda\leq
\lambda_{0}$. Moreover, $\displaystyle$%
\[
\lim_{\lambda\rightarrow0^{+}}\Vert u_{\lambda}\Vert_{L^{\infty}\left(
\Omega\right)  }=0.
\]

\end{theorem}

\textit{Proof}. Let $\delta>0$ and $w$ be as in the statement of the theorem. Define%

\[
\sigma:=\frac{1}{1+p}\in\left(  0,1\right)  .
\]
Note that $1-\sigma p=\sigma$. Let $v:=\lambda^{\sigma}w.$ So $v$ is a
subsolution of $(P_{\lambda})$ if and only if
\[
-\Delta v=\lambda^{\sigma}(m(x)-\delta)w^{-p}\leq\lambda m(x)\left(
f(\lambda^{\sigma}w)-k\right)  ,
\]
i.e.
\[
(m(x)-\delta)w^{-p}\leq\lambda^{1-\sigma}m(x)\left(  f(\lambda^{\sigma
}w)-k\right)  .
\]
Multiplying by $w^{p}$ and using that $1-\sigma=\sigma p$, we get the
equivalent condition
\[
\delta\geq m(x)\left(  1-f(\lambda^{\sigma}w)\left(  \lambda^{\sigma}w\right)
^{p}\right)  +km(x)\left(  \lambda^{\sigma}w\right)  ^{p}.
\]
Now, since $m$ and $w$ are bounded and $\displaystyle \lim_{s\rightarrow0^{+}%
}f(s)s^{p}=1$, we see that for $\lambda$ sufficiently small the inequality
above holds, and thus $v$ is a subsolution of $(P_{\lambda})$.

On the other hand, let $\varphi:=\mathcal{S}\left(  \left\vert m\right\vert
\right)  $ and $t>0$. Then,
\[
-\Delta\left(  t(\lambda\varphi)^{\sigma}\right)  =t\lambda^{\sigma}\left(
-\sigma\varphi^{\sigma-1}\Delta\varphi-\sigma\left(  \sigma-1\right)
\varphi^{\sigma-2}\left\vert \nabla\varphi\right\vert ^{2}\right)  \geq
t\sigma\lambda^{\sigma}\varphi^{\sigma-1}|m(x)|,
\]
so that $t(\lambda\varphi)^{\sigma}$ is a supersolution of $(P_{\lambda})$
whenever
\begin{equation}
t\sigma(\lambda\varphi)^{\sigma-1}\geq f(t(\lambda\varphi)^{\sigma}%
)\quad\text{and}\quad t\sigma(\lambda\varphi)^{\sigma-1}\geq k. \label{sup}%
\end{equation}
The first inequality in (\ref{sup}) holds if
\[
t^{p+1}\sigma=t^{p+1}\sigma(\lambda\varphi)^{\sigma(p+1)-1}\geq f(t(\lambda
\varphi)^{\sigma})\left(  t(\lambda\varphi)^{\sigma}\right)  ^{p}%
\qquad\text{in }\Omega.
\]
Recalling that $\lim_{s\rightarrow0^{+}}f(s)s^{p}=1$, we obtain some $s_{1}>0$
such that
\[
f(s)s^{p}<2\quad\text{for } 0<s\leq s_{1}.
\]
Consequently, if $t(\lambda\varphi)^{\sigma}\leq s_{1}$ and $t^{p+1}\sigma>2$
then
\[
f(t(\lambda\varphi)^{\sigma})\left(  t(\lambda\varphi)^{\sigma}\right)
^{p}<t^{p+1}\sigma.
\]
In other words, if we choose $t>0$ such that $t^{p+1}\sigma>2$ and
\[
\lambda_{0}\leq\left(  \frac{s_{1}}{t}\right)  ^{\frac{1}{\sigma}}\Vert
\varphi\Vert_{\infty}^{-1},
\]
then the first inequality in \eqref{sup} holds for $\lambda\leq\lambda_{0}$.
Moreover, if in addition%
\[
\lambda_{0}\leq\left(  \frac{t\sigma}{k}\right)  ^{\frac{1}{1-\sigma}}%
\Vert\varphi\Vert_{\infty}^{-1},
\]
then the second inequality in \eqref{sup} also is true for $\lambda\leq
\lambda_{0}$.

Finally, since $\varphi\in\mathcal{P}^{\circ}$, we may take $t$ in such a way
that $w\leq t\varphi^{\sigma}$. Therefore, \cite[Theorem 4.1]{loc} provides us
a solution $u_{\lambda}\in\mathcal{C}^{1}(\Omega)\cap\mathcal{C}%
(\overline{\Omega})$ of $(P_{\lambda})$ such that $\lambda^{\sigma}w\leq
u_{\lambda}\leq t\lambda^{\sigma}\varphi^{\sigma}$ for $0<\lambda\leq
\lambda_{0}$. In particular, we have that $\Vert u_{\lambda}\Vert_{\infty
}\rightarrow0$ as $\lambda\rightarrow0^{+}$. The proof is now complete.
$\blacksquare$

\end{document}